\documentclass[12pt]{article}
\usepackage{amssymb,latexsym}

\evensidemargin\oddsidemargin
\newtheorem{theorem}{Theorem}
\newtheorem{corollary}{Corollary}

\begin{document}

\bigskip

\begin{center}
{\LARGE Hilbert Norms For Graded Algebras}

\bigskip

{\large Joachim Kupsch\footnote{%
e-mail: kupsch@physik.uni-kl.de}}

{\large Fachbereich Physik der Universit\"at Kaiserslautern\\D-67663
Kaiserslautern, Germany}

\bigskip

{\large Oleg G. Smolyanov\footnote{%
e-mail: smolyan@mail.ru}}

{\large Faculty of Mechanics and Mathematics,\\ Moscow State
University, 119899 Moscow, Russia}

\medskip

\end{center}

\begin{abstract}
This paper presents a solution to a problem from superanalysis about the
existence of Hilbert-Banach superalgebras. Two main results are derived:\\1)
There exist Hilbert norms on some graded algebras (infinite-dimensional
superalgebras included) with respect to which the multiplication is
continuous.\\2) Such norms cannot be chosen to be submultiplicative and
equal to one on the unit of the algebra.\\\\AMS classification: 16W50,
46C05, 46H25
\end{abstract}

\section{Introduction}

The type of norms investigated in this article are generalizations of norms
used for the symmetric tensor algebra in the white noise analysis \cite
{HKPS:1993}\cite{Kon/Streit:1993} or in the Malliavin calculus \cite
{Watanabe:1984}. But now more general algebras are included, especially the
algebra of antisymmetric tensors (Grassmann algebra) and $\mathbb{Z}_2$%
-graded algebras (superalgebras) related to supersymmetry and to quantum
probability \cite{Meyer:1993}.

A locally convex commutative superalgebra is a $\mathbb{Z}_2$-graded locally
convex space $\mathcal{E}=\mathcal{E}_0\oplus \mathcal{E}_1$ equipped with
an associative continuous multiplication having the following property: for
any $a,b\in \mathcal{E}_0\cup \mathcal{E}_1,\,ab\neq 0$ the product
satisfies $ab=(-1)^{p(a)p(b)}ba$ with the parity function $p$, which is
defined on $\left( \mathcal{E}_0\cup \mathcal{E}_1\right) \setminus \left\{
0\right\} $ with $p\left( \mathcal{E}_0\setminus \left\{ 0\right\} \right)
=0,\,p\left( \mathcal{E}_1\setminus \left\{ 0\right\} \right) =1$, and $%
p(ab)=\left| p(a)-p(b)\right| $. Typical examples are Grassmann algebras
with finite or countable sets of generators. In superanalysis one considers
modules over (commutative) superalgebras \cite{Rogers:1980}\cite
{Jadczyk/Pilch:1981}\cite{DeWitt:1984}\cite{VV:1985b}\cite{Rogers:1986}\cite
{Berezin:1987}\cite{Smolyanov/Shavgulidze:1988}\cite{Khrennikov:1988}.%
\footnote{%
In the pioneering works of Martin \cite{Martin:1959} and of Berezin \cite
{Berezin:1966} the Grassmann algebra itself has been used instead of these
modules.} It is quite easy to define an infinite-dimensional Grassmann
algebra with a non-Hilbertian norm \cite{Rogers:1980}. But for a long time
it was unknown whether the topology of a locally convex superalgebra -
including the Grassmann algebra - can be defined with a Hilbert norm, and
moreover, whether this norm can be chosen to be simultaneously
submultilplicative and equal to one at the unit of the algebra. The paper
gives a complete solution to these problems. Our theorems imply a positive
answer to the first question and a negative answer to the second question.

\section{General considerations\label{gen}}

Let $\mathcal{A}$ be an algebra over the field $\mathbb{K}=\mathbb{R}$ or $%
\mathbb{C}$ with unit $e_0.$ The product is denoted by $a,b\in \mathcal{A}%
\rightarrow ab\in \mathcal{A}.$ We assume that $\mathcal{A}$ is provided
with a positive definite inner product $a,b\in \mathcal{A}\rightarrow \left(
a\mid b\right) \in \mathbb{K}.$ The corresponding Hilbert norm $\left\|
a\right\| =\sqrt{\left( a\mid a\right) }\geq 0$ is normalized at the unit $%
\left\| e_0\right\| =1.$ We are interested in such norms which allow a
uniform estimate for the product of the algebra 
\begin{equation}
\left\| ab\right\| \leq \gamma \left\| a\right\| \left\| b\right\|
\label{a.1}
\end{equation}
with a constant $\gamma \geq 1$. In this section we prove under rather
general conditions that this constant has the lower limit $\gamma \geq \sqrt{%
\frac 43}.$

\begin{theorem}
\label{Nogo1}Let $\mathcal{A}$ be an algebra over the field $\mathbb{K}=%
\mathbb{R}$ or $\mathbb{C}$ with dimension $\dim \mathcal{A}\geq 2.$ If this
algebra satisfies the properties\\ i)$\mathcal{A}$ is provided with a
Hilbert inner product $\left( .\mid .\right) $ normalized at the unit $e_0$, 
$\left\| e_0\right\| ^2=\left( e_0\mid e_0\right) =1$,\\ ii)there exists at
least one element $f\in \mathcal{A},f\neq 0$, such that $e_0,f$ and $f^2=ff$
satisfy $\left( e_0\mid f\right) =\left( f\mid f^2\right) =0$ and $\left(
e_0\mid f^2\right) \geq 0$, \\ then the norm estimate $\left\| ab\right\|
\leq \gamma \left\| a\right\| \left\| b\right\| $ is not valid for some $%
a,b\in \mathcal{A},$ if $\gamma <\sqrt{\frac 43}.$
\end{theorem}

\noindent \textbf{Proof }Since\textbf{\ }$f\neq 0$ we can normalize this
element and assume $\left\| f\right\| =1$. Take $a=e_0+\lambda f$ with $%
\lambda \in \mathbb{R}$. Then $a^2=e_0+2\lambda f+\lambda ^2f^2$ and \\$%
\left\| a^2\right\| ^2=1+2\lambda ^2\left( e_0\mid f^2\right) +4\lambda
^2+\lambda ^4\left\| f^2\right\| ^2\geq 1+4\lambda ^2.$ On the other hand $%
\left\| a\right\| ^2=1+\lambda ^2,$ and $\left\| a^2\right\| ^2\leq \gamma
^2\left\| a\right\| ^4$ implies $1+4\lambda ^2\leq \gamma ^2(1+\lambda
^2)^2. $ But this inequality is true for all $\lambda \geq 0$ only if $%
\gamma ^2\geq \sup_{\lambda \geq 0}(1+4\lambda ^2)(1+\lambda ^2)^{-2}=\frac
43$. \hfill $\Box $ \smallskip

This Theorem obviously applies to the tensor algebra $\mathcal{T}=\oplus
_{n=0}^\infty \mathcal{T}_n,$ where $\mathcal{T}_n$ is the subspace of
tensors of degree $n,$ and the norm is defined in the standard way as 
\begin{equation}
\left\| f\right\| ^2=\sum_{n=0}^\infty w_n\left\| f_n\right\| _n^2%
\mbox{  if
}\,f=\sum_{n=0}^\infty f_n,\;\,f_n\in \mathcal{T}_n  \label{a.2}
\end{equation}
with arbitrary positive weights $w_n>0,n\in \mathbb{N}$ and $w_0=1.$ In that
case we can simply choose an element $f\in \mathcal{T}_1,f\neq 0,$ to
satisfy the assumptions with $\left( e_0\mid f\otimes f\right) =0$.

Theorem \ref{Nogo1} can also be applied to a large class of algebras $%
\mathcal{A}$ which can be derived from the tensor algebra $\mathcal{T}$ by
the following modifications of the product.

\begin{enumerate}
\item  The product is generated by $f,g\in \mathcal{A}_1=\mathcal{T}%
_1\rightarrow f\circ g:=f\otimes g+(-1)^\chi g\otimes f$ where $\chi =0,1$ $%
\mathrm{mod}\;2$ is a parity factor.

\item  The product is generated by $f,g\in \mathcal{A}_1=\mathcal{T}%
_1\rightarrow \\f\circ g:=f\otimes g+(-1)^\chi g\otimes f+\omega (f,g)e_0.$
Here $\chi $ is again a parity factor and $\omega (.,.):\mathcal{A}_1\times 
\mathcal{A}_1\rightarrow \mathbb{K}$ is a bilinear pairing.
\end{enumerate}

The first class of algebras includes the algebra of symmetric tensors, the
algebra of antisymmetric tensors (Grassmann algebra), and tensor products of
these algebras, including the $\mathbb{Z}_2$-graded algebras (superalgebras)
used in quantum field theory. The assumptions of the Theorem\textbf{\ }\ref
{Nogo1} are satisfied for any non-vanishing element $f\in \mathcal{A}_1=%
\mathcal{T}_1.$

The second class includes the Clifford product, the (symmetric) Wiener
product, the antisymmetric Wiener product (with antisymmetric $\omega $) and
Le Jan's supersymmetric Wiener-Grassmann product \cite{LeJan:1988}\cite
{Kupsch:1990}\cite{Meyer:1993}. In these cases the assumptions of Theorem 
\ref{Nogo1} are satisfied if there exists a non-vanishing $f\in \mathcal{A}%
_1 $ with $\omega (f,f)\geq 0.$ Such a vector can always be found

if the algebra is complex, or

if the algebra is real and $\omega $ is not negative definite.

The last constraint is satisfied for the symmetric Wiener product on real
spaces, and for the real Clifford system in quantum field theory \cite
{BSZ:1992}. In both cases the form $\omega $ is positive definite.

Moreover Theorem \ref{Nogo1} is obviously true for any unital algebra $%
\mathcal{A}$, which has a nilpotent element $f$ that is orthogonal to the
unit element. If we only know that $\mathcal{A}$ has at least one nilpotent
element, we can derive the weaker

\begin{corollary}
Let $\mathcal{A}$ be an algebra which satisfies condition i) of Theorem \ref
{Nogo1}. If this algebra has a nilpotent element $f$, then the norm estimate 
$\left\| ab\right\| \leq \left\| a\right\| \left\| b\right\| $ is not valid
for some $a,b\in \mathcal{A}.$
\end{corollary}

\textbf{Proof} We assume again $\left\| f\right\| =1$. Then $a=e_0+\lambda f$
with $\lambda \in \mathbb{R}$ and $a^2=\left( e_0+\lambda f\right)
^2=e_0+2\lambda f$ have the norms $\left\| a\right\| ^2=1+2\lambda \mathrm{Re%
}\left( e_0,f\right) +\lambda ^2$ and $\left\| a^2\right\| ^2=1+4\lambda 
\mathrm{Re}\left( e_0,f\right) +4\lambda ^2.$ If $\mathrm{Re}\left(
e_0,f\right) =0$ we can apply the arguments given in the proof for Theorem 
\ref{Nogo1}. If $\mathrm{Re}\left( e_0,f\right) =\gamma \neq 0,$ then we
chose $\lambda =-2\gamma ,$ and $\left\| a^2\right\| ^2=1+8\gamma ^2\leq
1=\left\| a\right\| ^4$ is a contradiction.\hfill $\Box $ \smallskip 

\section{Norm estimates for $\mathbb{Z}$-graded algebras\label{NormEst}}

In this section we present Hilbert norm estimates for rather general $%
\mathbb{Z} $-graded algebras $\mathcal{A}$ over the field $\mathbb{K}=%
\mathbb{R}$ or $\mathbb{C}$.. We assume the following structure of 
$\mathcal{%
A}$.

\begin{enumerate}
\item  The algebra is the direct sum $\mathcal{A}=\oplus _{n=0}^\infty 
\mathcal{A}_n$ of orthogonal spaces $\mathcal{A}_n.$ Thereby $\mathcal{A}_0$
is the one dimensional space $\mathbb{K}$ spanned by the unit $e_0$ of the
algebra. The product $a\circ b$ maps $\mathcal{A}_p\times \mathcal{A}_q$
into $\mathcal{A}_{p+q}$ for all $p,q\in \left\{ 0,1,...\right\} .$

\item  The spaces $\mathcal{A}_n$ are provided with Hilbert norms $\left\|
.\right\| _n,n=0,1,....$The unit has norm $\left\| e_0\right\| _0=1.$ The
product of two homogeneous elements $a_p\in \mathcal{A}_p$ and $b_q\in 
\mathcal{A}_q$ satisfies 
\begin{equation}
\left\| a_p\circ b_q\right\| _{p+q}\leq \left\| a_p\right\| _p\left\|
b_q\right\| _q  \label{n.1}
\end{equation}
if $a_p\in \mathcal{A}_p$ and $b_q\in \mathcal{A}_q.$

\item  The algebra is provided with a family of Hilbert norms 
\begin{equation}
\left\| a\right\| _{(\sigma )}^2=\sum_{n=0}^\infty w_n(\sigma )\left\|
a_n\right\| _n^2\mbox{  if }\,a=\sum_{n=0}^\infty a_n,\;a_n\in \mathcal{A}_n
\label{n.2}
\end{equation}
with $\sigma \in \mathbb{R}.$ The factors $w_n(\sigma ),n=0,1,...,$ are
positive weights with the normalization $w_0(\sigma )=1$ for all $\sigma \in 
\mathbb{R}.$ The weights satisfy the inequalities $w_n(\sigma )\leq w_n(\tau
)$ for all $n\in \mathbb{N}$ if $\sigma \leq \tau .$
\end{enumerate}

\noindent An immediate consequence of these assumptions is $\left\|
a\right\| _{(\sigma )}\leq \left\| a\right\| _{(\tau )}$ for all $a\in 
\mathcal{A}$ if $\sigma \leq \tau .$ A simple example of such an algebra $%
\mathcal{A}$ is the tensor algebra $\mathcal{T}$. Its standard norm
satisfies (\ref{n.1}) with weights $w_n=1$ for all $n=0,1...$. More
interesting examples are the algebras of symmetric tensors or of
antisymmetric tensors. With the notation $f\circ g$ for both the symmetric
and the antisymmetric tensor product the estimate (\ref{n.1}) is satisfied
by the norms 
\begin{equation}
\left\| f_1\circ f_2\circ ...\circ f_n\right\| _n^2=\left\{ 
\begin{array}{c}
(n!)^{-1}\mathrm{per}(f_\mu \mid f_\nu )\mbox{  for symmetric tensors,}\, \\ 
(n!)^{-1}\mathrm{\det }\left( f_\mu \mid f_\nu \right) 
\mbox{  for
antisymmetric tensors,}\,
\end{array}
\right.  \label{n.1a}
\end{equation}
but it is violated if the factor $(n!)^{-1}$ is omitted. The standard norm%
\footnote{%
The ``standard'' inner product of the symmetric/antisymmetric tensor algebra
is characterized by the following property. Let $\mathcal{F}_i,i=1,2,$ be
two orthogonal subspaces of the space $\mathcal{A}_1.$ Denote by $\mathcal{A}%
(\mathcal{F}_i)$ the subalgebra generated by elements $f\in \mathcal{F}_i$.
Then $\left( a_1\circ a_2\mid b_1\circ b_2\right) =\left( a_1\mid b_1\right)
\left( a_2\mid b_2\right) $ holds for all $a_i\in \mathcal{A}(\mathcal{F}%
_i),i=1,2.$} is defined without the factor $(n!)^{-1}.$ In the notations
used here it corresponds therefore to a norm (\ref{n.2}) with a weight
function $w_n=n!$.

\begin{theorem}
\label{Inequality}If there exists a constant $\delta (\sigma ,\tau ,\rho )>0$
such that the weight functions satisfy the inequalities 
\begin{equation}
(p+q-1)w_{p+q}(\rho )\leq \delta (\sigma ,\tau ;\rho )w_p(\sigma )w_q(\tau )%
\mbox{  if }\,p,q\geq 1  \label{n.4}
\end{equation}
for values of $\sigma ,\tau $ and $\rho $ with $\sigma \leq \rho $ and $\tau
\leq \rho ,$ then the product of $\mathcal{A}$ is estimated by 
\begin{equation}
\left\| a\circ b\right\| _{(\rho )}\leq \gamma \cdot \left\| a\right\|
_{(\sigma )}\left\| b\right\| _{(\tau )}  \label{n.5}
\end{equation}
where the constant $\gamma $ is $\gamma =\sqrt{3}\max (1,\delta (\sigma
,\tau ,\rho ))$.
\end{theorem}

\noindent \textbf{Proof} For $a=a_0+a_{+}$ and $b=b_0+b_{+}$ with $%
a_0,b_0\in \mathcal{A}_0=\mathbb{K}$ and $a_{+}=\sum_{n=1}^\infty a_n$, $%
b_{+}=\sum_{n=1}^\infty b_n$ with $a_n,b_n\in \mathcal{A}_n,n\in \mathbb{N}$
the norm of $a\circ b$ is calculated by 
\[
\begin{array}{l}
\left\| a\circ b\right\| _{(\rho )}^2=\left\|
a_0b_0+a_0b_{+}+a_{+}b_0+a_{+}\circ b_{+}\right\| _{(\rho )}^2 \\ 
\leq \left| a_0b_0\right| ^2+3\left( \left| a_0\right| ^2\left\|
b_{+}\right\| _{(\rho )}^2+\left\| a_{+}\right\| _{(\rho )}^2\left|
b_0\right| ^2+\left\| a_{+}\circ b_{+}\right\| _{(\rho )}^2\right) \\ 
\leq \left| a_0b_0\right| ^2+3\left( \left| a_0\right| ^2\left\|
b_{+}\right\| _{(\rho )}^2+\left\| a_{+}\right\| _{(\rho )}^2\left|
b_0\right| ^2+\sum_{n\geq 1}w_n(\rho )\left\| \sum_{p+q=n}^{^{\prime
}}a_p\circ b_q\right\| _n^2\right)
\end{array}
\]
The symbol $\sum^{^{\prime }}$ means summation with the constraint $p\geq
1,q\geq 1$. The sum \\$\sum_{p+q=n,p\geq 1,q\geq
1}...=\sum_{p+q=n}^{^{\prime }}...$ has $n-1$ terms, hence \\$\left\|
\sum_{p+q=n}^{^{\prime }}a_p\circ b_q\right\| _n^2\leq
(n-1)\sum_{p+q=n}^{^{\prime }}\left\| a_p\circ b_q\right\| _n^2\stackrel{(%
\ref{n.1})}{\leq }(n-1)\sum_{p+q=n}^{^{\prime }}\left\| a_p\right\|
_p^2\left\| b_q\right\| _q^2$. \\If $w_n(\rho )$ is chosen such that (\ref
{n.4}) is satisfied we obtain \\$\sum_{n\geq 1}w_n(\rho )\left\|
\sum_{p+q=n}^{^{\prime }}a_p\circ b_q\right\| _n^2\leq \delta \cdot \left(
\sum_{p\geq 1}w_p(\sigma )\left\| a_p\right\| _p^2\right) \cdot \left(
\sum_{q\geq 1}w_q(\tau )\left\| b_q\right\| _q^2\right) $ \\$\leq \delta
\left\| a_{+}\right\| _{(\sigma )}^2\left\| b_{+}\right\| _{(\tau )}^2.$ For 
$\rho \leq \sigma ,\tau $ we have in addition the inequalities \\$\left\|
a_{+}\right\| _{(\rho )}^2\leq \left\| a_{+}\right\| _{(\sigma )}^2$ and $%
\left\| b_{+}\right\| _{(\rho )}^2\leq \left\| b_{+}\right\| _{(\tau )}^2$
such that finally 
\[
\begin{array}{ll}
\left\| a\circ b\right\| _{(\rho )}^2 & \leq \left| a_0b_0\right| ^2+3\left(
\left| a_0\right| ^2\left\| b_{+}\right\| _{(\tau )}^2+\left\| a_{+}\right\|
_{(\sigma )}^2\left| b_0\right| ^2+\delta \left\| a_{+}\right\| _{(\sigma
)}^2\left\| b_{+}\right\| _{(\tau )}^2\right) \\ 
& \leq 3\gamma \left\| a\right\| _{(\sigma )}^2\left\| b\right\| _{(\tau
)}^2.
\end{array}
\]
where $\gamma $ is $\gamma =\max (1,\delta ).$ \hfill $\Box $ \smallskip

As the first application of Theorem \ref{Inequality} we derive norms with
respect to which the product of the algebra is continuous. In that case the
inequality (\ref{n.4}) has to be satisfied for identical weights $w_p(\sigma
)=w_p(\tau )=w_p(\rho )=w_p,\;p\geq 1.$ If we fix $q=1$ then (\ref{n.4})
implies $p\cdot w_{p+1}\leq \delta \cdot w_p\cdot w_1$ for $p\in 
\mathbb{N}.$
As a consequence we obtain $w_p\leq \delta ^{p-1}\left( (p-1)!\right)
^{-1}w_1,p\geq 1.$ The slowest decrease of the weights which might be
possible according to our estimates is therefore $w_p\sim \left(
(p-1)!\right) ^{-1}.$ The proof that such a solution actually exists follows
from the simple estimate $\left( 
\begin{array}{c}
m+n \\ 
m
\end{array}
\right) =\frac{(m+n)!}{m!n!}\geq 1$ if $m,n\geq 0.$ Hence $(p+q-1)\frac
1{(p+q-1)!}=\frac 1{(p+q-2)!}\leq \frac 1{(p-1)!}\frac 1{(q-1)!}$ is valid
for all $p,q\geq 1.$ Since 
\begin{equation}
2^{m+n}\geq \left( 
\begin{array}{c}
m+n \\ 
m
\end{array}
\right) =\frac{(m+n)!}{m!n!}\geq m+n\mbox{  if }\,m,n\geq 1,  \label{n.6}
\end{equation}
also $(p+q-1)\frac 1{(p+q)!}<\frac 1{(p+q-1)!}\leq \frac 1{p!}\frac 1{q!}$
follows for all $p,q\geq 1.$ We have therefore derived

\begin{corollary}
\label{Estim1}If the norm is defined with the weights $w_0=1,\;w_n=\left(
(n-1)!\right) ^{-1},$ $n\geq 1,$ or with $w_0=1,w_n=\left( n!\right)
^{-1},\;n\geq 1,$ the product of the algebra is continuous with the uniform
norm estimate 
\begin{equation}
\left\| a\circ b\right\| \leq \sqrt{3}\left\| a\right\| \left\| b\right\| .
\label{n.7}
\end{equation}
\end{corollary}

As a more general class of norms we choose weights 
\begin{equation}
w_0=1,\;w_n(\sigma ,\rho ,s)=(n!)^\sigma 2^{\rho n}(1+n)^s\mbox{  if }%
\,n\geq 1,  \label{n.8}
\end{equation}
with real parameters $\sigma ,\rho ,s.$ These weights satisfy the
inequalities \\$w_n(\sigma _1,\rho _1,s_1)\leq w_n(\sigma _2,\rho _2,s_2)$
if $\sigma _1\leq \sigma _2,\rho _1\leq \rho _2,s_1\leq s_2$. We denote by $%
\left\| a\right\| _{(\sigma ,\rho ,s)}$ the norm (\ref{n.2}) defined with
the weights $w_n(\sigma ,\rho ,s).$ The estimate (\ref{n.6}) and the bounds $%
\frac{(m+n)!}{m!n!}\geq \frac{(2m)!}{(m!)^2}\geq const\cdot 2^{2m}m^{-\frac
12}$ if $n\geq m\geq 1$ and $1\leq \frac{(1+m)(1+n)}{1+m+n}\leq 1+\mathrm{%
\min }(m,n)$ yield inequalities of the type (\ref{n.4}) also for these
norms. We obtain 
\begin{equation}
(p+q-1)w_{p+q}(\sigma ,\rho ,s)\leq \delta w_p(\sigma ^{\prime },\rho
^{\prime },s^{\prime })w_q(\sigma ^{\prime },\rho ^{\prime },s^{\prime })%
\mbox{  if }\,p,q\geq 1  \label{n.9}
\end{equation}
with a constant $\delta \geq 1$ if $\sigma =\sigma ^{\prime }=-1$ with $\rho
=\rho ^{\prime }\in \mathbb{R}$ and $s=s^{\prime }\leq 0,$ or if $\sigma
=\sigma ^{\prime }<-1$ with $\rho =\rho ^{\prime }\in \mathbb{R}$ and $%
s=s^{\prime }\in \mathbb{R}$.

The generalizations of (\ref{n.7}) are therefore 
\begin{equation}
\left\| a\circ b\right\| _{(-1,\rho ,s)}\leq \sqrt{3}\left\| a\right\|
_{(-1,\rho ,s)}\cdot \left\| b\right\| _{(-1,\rho ,s)}\mbox{  if }\,\rho \in 
\mathbb{R},s\leq 0,  \label{n.10a}
\end{equation}
and 
\begin{equation}
\left\| a\circ b\right\| _{(\sigma ,\rho ,s)}\leq \gamma \left\| a\right\|
_{(\sigma ,\rho ,s)}\cdot \left\| b\right\| _{(\sigma ,\rho ,s)}\mbox{  if }%
\,\sigma <-1,\rho \in \mathbb{R},s\in \mathbb{R.}  \label{n.10b}
\end{equation}
Here $\gamma $ takes some value $\gamma \geq \sqrt{3}$ depending on the
choice of the parameters $\sigma $ and $s.$

Moreover, the inequalities (\ref{n.9}) are valid for $(\sigma ,\rho ,s)\neq
(\sigma ^{\prime },\rho ^{\prime },s^{\prime })$ if $\sigma <\sigma ^{\prime
}$ or if $\sigma =\sigma ^{\prime }$ and $\rho <\rho ^{\prime }$. The
corresponding estimates for the norms are 
\begin{equation}
\left\| a\circ b\right\| _{(\sigma ,\rho ,s)}\leq \gamma \left\| a\right\|
_{(\sigma ^{\prime },\rho ^{\prime },s^{\prime })}\cdot \left\| b\right\|
_{(\sigma ^{\prime },\rho ^{\prime },s^{\prime })}\mbox{  if }\,\sigma
<\sigma ^{\prime }\mbox{  for all }\,\rho ,\rho ^{\prime },s,s^{\prime }\in 
\mathbb{R,}  \label{n.11a}
\end{equation}
and 
\begin{equation}
\left\| a\circ b\right\| _{(\sigma ,\rho ,s)}\leq \gamma \left\| a\right\|
_{(\sigma ,\rho ^{\prime },s^{\prime })}\cdot \left\| b\right\| _{(\sigma
,\rho ^{\prime },s^{\prime })}\mbox{  if }\,\rho <\rho ^{\prime }%
\mbox{  for all 
}\,\sigma ,s,s^{\prime }\in \mathbb{R.}  \label{n.11b}
\end{equation}
The value of $\gamma \geq \sqrt{3}$ depends on the choice of the parameters.

For the tensor algebra and for algebras of symmetrized tensors\footnote{%
This class of algebras does not only include the algebra of symmetric
tensors and the algebra of antisymmetric tensors (Grassmann algebra), but
also the $\mathbb{Z}_2$-graded algebras (superalgebras) used in
supersymmetric quantum field theory.} the Hilbert space $\mathcal{A}_1=%
\mathcal{H}$ generates the whole algebra. Given a (self-adjoint/positive)
operator $A$ on $\mathcal{H}$, the mapping $\Gamma (A)e_0=e_0$ and $\Gamma
(A)\left( f_1\circ f_2\circ ...\circ f_n\right) :=(Af_1)\circ (Af_2)\circ
...\circ (Af_n)$ for $f_\mu \in \mathcal{H},\mu =1,...,n,$ and $n\in 
\mathbb{%
N}$, defines a unique (self-adjoint/positive) operator $\Gamma (A)$ on the
algebra $\mathcal{A}$, which satisfies the relation 
\begin{equation}
\Gamma (A)(a\circ b)=\left( \Gamma (A)a\right) \circ \left( \Gamma
(A)b\right) .  \label{n.12}
\end{equation}
The norms (\ref{n.2}) with the weights (\ref{n.8}) are then easily
generalized to 
\begin{equation}
\left\| a\right\| _{(\sigma ,\rho ,s)}^2=\sum_{n=0}^\infty (n!)^\sigma
\left\| \left( \Gamma (A)\right) ^\rho a_n\right\| _n^2(1+n)^s\mbox{  if }%
\,a=\sum_{n=0}^\infty a_n,\;a_n\in \mathcal{A}_n.  \label{n.13}
\end{equation}
If $A$ is an invertible positive operator with lower bound $A\geq 2\cdot id$%
, then $\Gamma (A)$ satisfies $\left\| \left( \Gamma (A)\right) ^{-\rho
}a\right\| _n\leq 2^{-n\rho }\left\| a\right\| _n$ for $a\in \mathcal{A}_n$
if $\rho \geq 0$. This bound and the relation (\ref{n.12}) imply that the
estimates (\ref{n.10a}),(\ref{n.10b}) and (\ref{n.11b}) are also valid for
the norms (\ref{n.13}), moreover (\ref{n.11a}) holds if $\rho \leq \rho
^{\prime }$.

If $A^{-1}$ is a Hilbert-Schmidt operator then a family of norms (\ref{n.13}%
) can be used to define a nuclear topology on the algebra $\mathcal{A}$. For
the symmetric tensor algebra that has been done in the white noise calculus
and in the Malliavin calculus, see e.g. \cite{Arai/Mitoma:1993} \cite
{Kon/Streit:1993} \cite{Watanabe:1984}. For the algebra of antisymmetric
tensors and for the superalgebras such nuclear topologies can be found in 
\cite{Kree:1978a} and in \cite{Haba/Kupsch:1995}. But the estimates of these
references are not strong enough to derive the results with a single Hilbert
norm as presented in Corollary \ref{Estim1} and in eqs. (\ref{n.10a}) and (%
\ref{n.10b}).

\begin{center}
\textbf{Acknowledgment}
\end{center}

\noindent A great part of this work was done during a stay of O. G.
Smolyanov at the University of Kaiserslautern. OGS would like to thank the
Deutsche Forschungsgemeinschaft (DFG) and the Russian Fund of Fundamental
Research for financial support.


\begin{thebibliography}{10}

\bibitem{Arai/Mitoma:1993}
A.~Arai and I.~Mitoma.
\newblock Comparison and nuclearity of spaces of differential forms on
  topological vector spaces.
\newblock {\em J. Funct. Anal.}, 111:278--294, 1993.

\bibitem{BSZ:1992}
J.~C. Baez, I.~E. Segal, and Z.~Zhou.
\newblock {\em {Introduction to Algebraic and Constructive Quantum Field
  Theory}}.
\newblock Princeton University Press, Princeton, 1992.

\bibitem{Berezin:1966}
F.~A. Berezin.
\newblock {\em The Method of Second Quantization}.
\newblock Academic Press, New York, 1966.

\bibitem{Berezin:1987}
F.~A. Berezin.
\newblock {\em {Introduction to Superanalysis}}.
\newblock Reidel, Dordrecht, 1987.

\bibitem{DeWitt:1984}
B.~DeWitt.
\newblock {\em Supermanifolds}.
\newblock CUP, Cambridge, 1984.

\bibitem{Haba/Kupsch:1995}
Z.~Haba and J.~Kupsch.
\newblock Supersymmetry in euclidean quantum field theory.
\newblock {\em Fortschr. Phys.}, 43:41--66, 1995.

\bibitem{HKPS:1993}
T.~Hida, H.-H. Kuo, J.~Potthoff, and L.~Streit.
\newblock {\em White Noise}.
\newblock Kluwer, Dordrecht, 1993.

\bibitem{Jadczyk/Pilch:1981}
A.~Jadczyk and K.~Pilch.
\newblock Superspaces and supersymmetries.
\newblock {\em Commun. Math. Phys.}, 78:373--390, 1981.

\bibitem{LeJan:1988}
Y.~Le Jan.
\newblock {On the Fock space representation of functionals of the occupation
  number field and their renormalization}.
\newblock {\em J. Funct. Anal.}, 80:88--108, 1988.

\bibitem{Khrennikov:1988}
A.~Yu. Khrennikov.
\newblock Functional superanalysis.
\newblock {\em Russian Math. Surveys}, 43:103--137, 1988.

\bibitem{Kon/Streit:1993}
Yu~G. Kondratev and L.~Streit.
\newblock Spaces of white noise distributions: constructions, descriptions,
  applications. {I}.
\newblock {\em Rep. Math. Phys.}, 33:341--366, 1993.

\bibitem{Kree:1978a}
P.~Kr{\'e}e.
\newblock M{\'e}thodes fonctionelles en analyse de dimension infinie et
  holomorphie anticommutative.
\newblock In {\em {S{\'e}minaire P. Lelong et H. Skoda (Analyse) Ann{\'e}e
  1976/77, Lect. Notes Math. 694}}, pages 134--171. Springer, Berlin, 1978.

\bibitem{Kupsch:1990}
J.~Kupsch.
\newblock A probabilistic formulation of bosonic and fermionic integration.
\newblock {\em Rev. Math. Phys.}, 2:457--477, 1990.

\bibitem{Martin:1959}
J.~L. Martin.
\newblock {Generalized classical dynamics, and the "classical analogue" of a
  Fermi oscillator}.
\newblock {\em Proc. Roy. Soc. (London)}, A251:536--542, 1959.

\bibitem{Meyer:1993}
P.~A. Meyer.
\newblock {\em Quantum Probability for Probabilists}.
\newblock Lect. Notes in Math. 1538. Springer, Berlin, 1993.

\bibitem{Rogers:1980}
A.~Rogers.
\newblock A global theory of supermanifolds.
\newblock {\em J. Math. Phys.}, 21:1352--1365, 1980.

\bibitem{Rogers:1986}
A.~Rogers.
\newblock Graded manifolds, supermanifolds and infinite-dimensional {Grassmann}
  algebras.
\newblock {\em Commun. Math. Phys.}, 105:375--384, 1986.

\bibitem{Smolyanov/Shavgulidze:1988}
O.~G. Smolyanov and E.~T. Shavgulidze.
\newblock {The Fourier transform and pseudodifferential operators in
  superanalysis}.
\newblock {\em Soviet Math. Dokl.}, 37:476--481, 1988.

\bibitem{VV:1985b}
V.~S. Vladimirov and I.~V. Volovich.
\newblock On the definition of the integral in superspace.
\newblock {\em Soviet Math. Dokl.}, 32:817--819, 1985.

\bibitem{Watanabe:1984}
S.~Watanabe.
\newblock {\em {Lectures on Stochastic Differential Equations and Malliavin
  Calculus}}.
\newblock Tata Institute, Bombay, and Springer-Verlag, Berlin, 1984.

\end{thebibliography}
\end{document}